\documentclass[11pt]{article}
\usepackage{fullpage}
\usepackage{amsfonts, amsmath, amssymb,latexsym,epsfig}
\usepackage{amsthm}
\usepackage{tikz}
\usetikzlibrary{graphs}
\usetikzlibrary{graphs.standard}
\usepackage{relsize}
\usepackage{verbatim}
\usepackage{enumerate}

\usepackage{tabularx}
\usepackage{adjustbox}
\usepackage[skip=10pt plus1pt, indent=20pt]{parskip}

\newtheorem{thm}{Theorem}[section]
\newtheorem{cor}[thm]{Corollary}

\newtheorem{defn}[thm]{Definition}

\providecommand{\keywords}[1]{\textbf{\textit{Keywords}} #1}

\title{Chromatic and achromatic numbers of unitary addition Cayley graphs}

\author{Keenan Calhoun\thanks{Department of Mathematics and Statistics, California State University Sacramento.
\texttt{kcalhoun@csus.edu}} \and 
Ye\c{s}\.im Dem\.iro\u{g}lu Karabulut\thanks{Department of Mathematics and Statistics, California State University Sacramento.   
		\texttt{demiroglu@csus.edu}}
  \and
  Vincent Pigno\thanks{Department of Mathematics and Statistics, California State University Sacramento.   
		\texttt{vincent.pigno@csus.edu}}
  \and
Craig Timmons\thanks{Department of Mathematics and Statistics, California State University Sacramento.   
		\texttt{craig.timmons@csus.edu}}}

\begin{document}
	
	\maketitle

	\begin{abstract}
		 Let $R$ be a ring.  The unitary addition Cayley graph of $R$, denoted $\mathcal{U}(R)$, is the graph with vertex $R$, and two distinct vertices $x$ and $y$ are adjacent if and only if $x+y$ is a unit.  We determine a formula for the clique number and chromatic number of such graphs when $R$ is a finite commutative ring with an odd number of elements.  This includes the special case when $R$ is $\mathbb{Z}_n$, the integers modulo $n$, where these parameters had been found under the assumption that $n$ is even, or $n$ is a power of an odd prime. 
   Additionally, we study the achromatic number 
   of $\mathcal{U}( \mathbb{Z}_n )$ in the case that $n$ is the product of two primes.  We prove that the achromatic number of $\mathcal{U} ( \mathbb{Z}_{3q})$ is equal to $\frac{3q+1}{2}$ when $q > 3$ is a prime.  We also prove a lower bound that applies when $n = pq$ where $p$ and $q$ are distinct odd primes.    
	\end{abstract}
	


\section{Introduction}

Throughout this paper, all rings are assumed to be finite rings with unity.  There are many different graphs that can be associated to a ring.  One may use combinatorial properties of the graph to study algebraic properties of the ring, or perhaps choose the ring so that the graph constructed has some desired property.  The recent monograph of Anderson, Asir, Badawi, and Chelvam \cite{Anderson} discusses several graphs obtained from rings such as zero-divisor graphs, total graphs, and Cayley graphs.  The subject of this work is unitary addition Cayley graphs and their clique, chromatic and achromatic numbers.  These graphs have been the subject of many papers \cite{Akhtar, Bahrami, Chin, Klotz, Liu, Momrit, Palanivel, Promsakon, Sinha}.  Their definition is as follows.  

Let $R$ be a ring and $R^*$ be the units of $R$.  
The \emph{unitary addition Cayley graph} of $R$, denoted by $\mathcal{U}(R) $, is the graph whose vertex set is $R$, and distinct vertices $x$ and $y$ are adjacent if and only if $x + y \in R^*$.  One of the most studied cases is when $R =  \mathbb{Z}_n$, the cyclic group of integers modulo $n$.
Basic properties of $\mathcal{U}(  \mathbb{Z}_n )$, such as the number of edges, were determined by Sinha, Garg, and Singh \cite{Sinha}.  In the case that $n$ is even, 
they proved $\mathcal{U} ( \mathbb{Z}_n )$
is a $\phi (n)$-regular bipartite graph.  
Here $\phi (n)$ denotes the \emph{Euler Phi-function} and is the number of integers in $\{1,2, \dots , n \}$ that are relatively prime with $n$.  For odd $n$, the degree of a vertex $x$ in 
$\mathcal{U} ( \mathbb{Z}_n )$ is 
\[  
d (x) = 	 \left\{
	  	\begin{array}{ll}
	  		\phi (n) & \mbox{if $n$ is odd and $\textup{gcd}(x,n) > 1$,} \\
	  		\phi (n) - 1  & \mbox{if $n$ is an odd and $\textup{gcd}(x,n) = 1$.} 
	  	\end{array} 
	  	\right. \]
In 2007, Klotz and Sander \cite{Klotz} proved a number of results about the graphs $\mathcal{U}( \mathbb{Z}_n)$.  Since their work, many graphs parameters such as the girth, Hamiltonicity, planarity, and the spectrum of $\mathcal{U}( \mathbb{Z}_n)$ have been studied.   
        Let us now focus on what is known regarding the clique number, chromatic number, and achromatic number of 
    $\mathcal{U} ( \mathbb{Z}_n )$.  
 Observe that $\chi ( \mathcal{U} ( \mathbb{Z}_n ) ) 
 =2$ since $\mathcal{U}( \mathbb{Z}_n )$ is bipartite when $n$ is even.  Sinha, Garg, and Singh \cite{Sinha} proved that for $p$ an odd prime and an odd integer $n$ having $m \geq 2$ distinct prime factors, 
 \begin{center}
     $\chi ( \mathcal{U} ( \mathbb{Z}_p ) )  = \dfrac{p +1}{2}$ ~~ and ~~ $\chi ( \mathcal{U} ( \mathbb{Z}_n )) \leq \dfrac{ \phi (n) }{2^m} + m$.
      \end{center}
   
These same results were obtained by Palanivel and Chithra
\cite{Palanivel} who also investigated $\mathcal{U} ( \mathbb{Z}_n )$. For a graph $G$, let $\alpha (G)$ and $\omega (G)$ denote the independence number and clique number of $G$, respectively.
Again, since $\mathcal{U} ( \mathbb{Z}_n )$ is bipartite when $n$ is even, we have $\omega ( \mathcal{U} ( \mathbb{Z}_n )) =2$.  
Palanivel and Chithra \cite{Palanivel} proved that for an odd prime $p$ and integer $k \geq 2$, 
\begin{center}
    $\omega ( \mathcal{U} ( \mathbb{Z}_{p^k}  )) =\dfrac{ \phi (p^k) }{2} + 1 $.
\end{center}
	  They also proved the lower bound 
	  	  \[
	  \omega ( \mathcal{U} ( \mathbb{Z}_n )) \geq 
	  \left\{ 
	  \begin{array}{ll} 
	  	3 & \mbox{if $p_1 = 3$},  \\
	  	\frac{p_1+1}{2} & \mbox{if $p_1 > 3$}
	  	\end{array}
  	\right.
	  \]
   where $n = p_1^{k_1} p_2^{k_2} \cdots p_m^{k_m} $ is the prime factorization of $n$ and $p_1 < p_2  < \dots < p_m$.

   By putting these results together, we get that  
 for $n = p^k$ where $p$ is an odd prime and $k \geq 1$ an integer, 
  \[
  \frac{ \phi (p^k) }{2} + 1 
  =
  \omega ( \mathcal{U} ( \mathbb{Z}_{p^k} ) ) \leq \chi (  \mathcal{U} ( \mathbb{Z}_{p^k} ) ) 
  \leq 
  \frac{ \phi ( p^k ) }{2} + 1.
  \]
  As the lower and upper bounds match, equality holds throughout.   
  The chromatic number of $\mathcal{U} ( \mathbb{Z}_{n})$ was studied further by Promsakon \cite{Promsakon} who also obtained
  some of the results of \cite{Palanivel}. 
  Thus, the value of $\chi( \mathcal{U} ( \mathbb{Z}_n ))$ has been determined when $n$ is even, and when $n$ is a power of an odd prime.  However, when $n$ has more than one prime factor exact formulas for the clique and chromatic number are not known.  Our first theorem gives such formulas for all odd $n$.  Furthermore, it applies to a large class of rings that includes $\mathbb{Z}_n$ as a special case.  Recall that a commutative ring with unity is a \emph{local ring} if it has a unique maximal ideal.  One of the classical results in ring theory is that every finite commutative ring is a product of local rings.       

\begin{thm}\label{theorem:general clique and chromatic}
    Let $R$ be a finite commutative ring with an odd number of elements that is the direct product of $m$ local rings, say 
    $R = R_1 \times  \cdots \times R_m$ where $M_i$ denotes the unique maximal ideal in $R_i$.  
    Then
    \[
   m +  \prod_{i=1}^m \left( \frac{ |R_i| - |M_i | } {2} \right)  = \omega ( \mathcal{U} ( R_1 \times \cdots \times R_m)) =
   \chi ( \mathcal{U} ( R_1 \times \cdots \times R_m)) .
    \]
\end{thm}

In the case of $\mathbb{Z}_n$ 
with $n$ an odd integer having prime factorization $n = p_1^{k^1} \cdots p_m^{k_m}$, we have the ring isomorphism $\mathbb{Z}_n \cong
\mathbb{Z}_{p_1^{k_1} } \times \cdots \times 
\mathbb{Z}_{p_m^{k_m} } $.  
Let us assume that $3 \leq p_1 < p_2 < \dots < p_m$.  
Each factor $\mathbb{Z}_{ {p_i}^{k_i} }$ is a local ring where the unique maximal ideal is $\langle p_i \rangle$.  
Therefore, $| \mathbb{Z}_{ {p_i}^{k_i} } |  - 
| \langle p_i \rangle | = p_{i}^{k_i} - p_i^{k_i - 1} = \phi ( p_i^{k_i} )$.  Using the multiplicative property of $\phi$, 
\[
m + \prod_{i=1}^{m} \left( \frac{ \phi ( p_i^{k_i} )  }{2} \right) = m + \frac{ \phi (n) }{2^m}.  
\]
Thus, by Theorem \ref{theorem:general clique and chromatic} we have the following Corollary.  

\begin{cor}\label{cor:1}
If $n$ is an odd positive integer with $m$ distinct prime factors, then 
\[
\omega ( \mathcal{U} ( \mathbb{Z}_n ) )  
= 
\chi ( \mathcal{U}  ( \mathbb{Z}_n ) ) 
= 
m + \frac{ \phi (n) }{2^m}. 
\]
\end{cor}

The previously best known bounds on the clique and chromatic number on $\mathcal{U}( \mathbb{Z}_n)$, as well as our contributions, are summarized in the following table.  \textbf{In rows 2, 3, 4, and 5}, $p, p_1, \dots , p_m$ are odd primes with $3 \leq p_1 < p_2 < \dots  < p_m$ and $k$ is a positive integer.    

\begin{center}
	\begin{tabular}{c|c|c|c}
		number of vertices & $\omega$ & $\chi$ & Reference \\ \hline
		$n$ where $n$ is even & $2 = \omega ( \mathcal{U}( \mathbb{Z}_n ) ) $ & $\chi ( \mathcal{U}( \mathbb{Z}_n ) ) = 2$ & $\mathcal{U} ( \mathbb{Z}_n)$ is bipartite \\
		$p$  & $ \frac{p+1}{2} = \omega ( \mathcal{U}( \mathbb{Z}_p ) ) $ 
		& $\chi ( \mathcal{U}( \mathbb{Z}_p ) ) = \frac{p+1}{2} $ & \cite{Promsakon, Sinha} \\
		$p^k$  & $\frac{\phi(p^k)}{2}+1 = \omega ( \mathcal{U}( \mathbb{Z}_{p^k} ) ) $ 
		& $\chi ( \mathcal{U}( \mathbb{Z}_{p^k} ) ) = \frac{ \phi( p^k )}{2}+1 $ & \cite{Palanivel, Promsakon}
		 \\
		$n=p_1^{k_1} \cdots p_m^{k_m}$ & $\frac{p_1+1}{2} \leq \omega ( \mathcal{U}( \mathbb{Z}_n ) ) $ 
		& $\chi ( \mathcal{U}( \mathbb{Z}_n ) ) \leq \frac{ \phi(n) }{2^m}+m $ & \cite{Palanivel, Promsakon, Sinha} \\ 
			$n=p_1^{k_1} \cdots p_m^{k_m}$ & $  \frac{ \phi(n) }{2^m}+m =\omega ( \mathcal{U}( \mathbb{Z}_n ) ) $ 
		& $\chi ( \mathcal{U}( \mathbb{Z}_n ) ) = \frac{ \phi(n) }{2^m}+m$ & Corollary \ref{cor:1} 
		\end{tabular}
        
        \textbf{Current Bounds on $\omega ( \mathcal{U} ( \mathbb{Z}_n )$ and $\chi ( \mathcal{U} ( \mathbb{Z}_n )$ and Corollary \ref{cor:1}}
\end{center}

Other unitary graphs, aside from $\mathcal{U}( \mathbb{Z}_n)$ have been investigated.  For example, Theorem \ref{theorem:general clique and chromatic} gives the clique and chromatic number of the unitary addition graphs of the Gaussian integers modulo $n$, where $n$ is odd.  These graphs were studied by Bahrami and Jahani-Nezhad (\cite{Bahrami}, Section 3.3).     

Originally motivated by a coloring problem on hypercubes, discussed in further detail in Section \ref{se:concluding remarks},
we also investigated the achromatic number of unitary addition Cayley graphs.  The \emph{achromatic number of a graph $G$}, denoted $\chi_a (G) $, is the maximum number of colors in a proper coloring where there is at least one edge between each pair of color classes.  
In 1967, Harary, Hedetniemi, and Prins introduced these types of colorings in connection with the chromatic number of a graph  \cite{Harary}.  
Achromatic colorings have been studied extensively since then \cite{Edwards, EdwardsSurvey, Harary2}.  In general, they are difficult to find.    
Not surprisingly, less is known about $\chi_a(  \mathcal{U} ( \mathbb{Z}_n ) )$ compared to
$\chi(  \mathcal{U} ( \mathbb{Z}_n ) )$.  
Momrit and Promsakon \cite{Momrit} proved that 
for any odd prime $p$ and integer $k \geq 1$, 
\[
\chi_a( \mathcal{U} ( \mathbb{Z}_{2^k} )) = 2 ~~~~ \mbox{and} 
~~~~
 \chi_a ( \mathcal{U} ( \mathbb{Z}_{p^k} )  ) = \frac{ \phi ( p^k ) }{2} + 1 
 ~~~~ \mbox{and} ~~  p \leq \chi_a ( \mathcal{U} ( \mathbb{Z}_{2p} ) ) \leq p + 1.
 \]
Thus, the achromatic number of
$\mathcal{U}( \mathbb{Z}_n )$ is known only when $n$ is a power of a prime.    
 In some cases, the chromatic number and the achromatic number of $\mathcal{U} ( \mathbb{Z}_n)$ are equal.  This occurs when $n$ consists of only one prime factor:  
    \[
    \chi ( \mathcal{U}  ( \mathbb{Z}_{2^k} ) ) = \chi_a ( \mathcal{U} ( \mathbb{Z}_{2^k} ) ) = 2 
    ~~~ \mbox{and} ~~~
    \chi (\mathcal{U} ( \mathbb{Z}_{p^k} ) ) = \chi_a( \mathcal{U} ( \mathbb{Z}_{p^k} ) ) = \frac{ \phi (p^k) }{2} + 1.
    \] 
    When $n$ is not a prime power, these values can differ substantially.  
    For example, while $\chi ( \mathcal{U} ( \mathbb{Z}_{2p} )) = 2$, we have $\chi_a ( \mathcal{U}  ( \mathbb{Z}_{2p} )) = p$.
In Section \ref{section:achromatic} we 
show that $\chi_a ( \mathcal{U}  ( \mathbb{Z}_{2p})) = p$ and prove the following lower bound that applies when $n$ is the product of two odd primes.     

\begin{thm}\label{thm:achromatic theorem 1}
If $q > p \geq 3$ are primes, then 
\[
\chi_a ( \mathcal{U} ( \mathbb{Z}_{pq} ) ) \geq  \frac{pq+1}{2}.
\]
\end{thm}

 We were able to prove that this lower bound is best possible when the smaller prime is 3.     

\begin{thm}\label{thm:achromatic theorem 2}
\noindent
\medskip
If $q > 3$ is prime, then $\chi_a ( \mathcal{U} ( \mathbb{Z}_{3q} ) ) =  \frac{3q+1}{2}$.  
\end{thm}

Theorem \ref{thm:achromatic theorem 2} determines the achromatic number for a non-trivial infinite family of Cayley sum graphs.  
Initially it was thought that the formula in Theorem \ref{thm:achromatic theorem 2} could hold more generally and we conjectured that the achromatic number of $\mathcal{U}( \mathbb{Z}_{pq} )$ is $\frac{pq+1}{2}$ for $q > p \geq 3$ prime.  
If this were true, then it would show that the bound from Theorem \ref{thm:achromatic theorem 1} is best possible. 
In addressing questions raised by the referee report(s), we discovered that this conjecture is false.
Using a computer search, we were able to find an achromatic coloring of $\mathcal{U}( \mathbb{Z}_{35} )$ with 19 colors, exceeding $\frac{5\cdot 7 + 1}{2}$.  An achromatic coloring of $\mathcal{U}( \mathbb{Z}_{55})$ with $30 > \frac{ 5\cdot 11 + 1}{2}$ colors was also found (the authors would like to thank the reviewer(s) for comments regarding our conjecture).  

   In general, it seems interesting to study $\chi_a ( \mathcal{U} ( R) )$ in the case that $R$ is a finite local ring.  Our investigation of $\chi_a ( \mathcal{U} ( \mathbb{Z}_n))$ was motivated by a problem of Ahlswede, Bezrukov, Blokhuis, Metsch, and  Moorhouse \cite{Ahlswede}.  This problem is discussed further in the Section \ref{se:concluding remarks}.



\section{Proof of Theorem \ref{theorem:general clique and chromatic}}\label{section:chromatic} 

\begin{proof}[Proof of Theorem \ref{theorem:general clique and chromatic}]
Let $R$ be a finite commutative ring with unity having an odd number of elements.  Assume that
$R$ is a product of local rings, say 
$R = R_1 \times  \cdots \times R_m$.  
Let $M_i$ be the unique maximal ideal in $R_i$.  
It is known that the units of $R_i$ are the elements of $R_i \backslash M_i$, and $\prod_{i=1}^m (R_i \backslash M_i)$ is the set of units in $R$.  Observe then that if $x = (x_1 , \dots , x_m)$ and $y = (y_1 , \dots , y_m )$ are distinct elements in $R$, then $x+y$ is a non-unit if and only if 
$x_i + y_i \in M_i$ for some $1 \leq i \leq m$.  
In terms of the graph 
$\mathcal{U} ( R_1 \times  \cdots \times R_m )$, a set of vertices $\mathcal{A}$ forms an independent set if for each $x \neq y \in \mathcal{A}$, the sum  $x+y$ has at least one coordinate in one of the maximal ideals. 
 We will use this fact frequently in our argument.  

First we prove the upper bound
\begin{equation}\label{eq:general ub}
 \chi ( \mathcal{U} ( R_1 \times \cdots \times R_m))
\leq 
m +  \prod_{i=1}^m \left( \frac{ |R_i| - |M_i | } {2} \right).
 \end{equation} 
Let $\tau = \frac{ |R^* | }{2^m }$.  
We will define a partition of $R$ into $m + \tau$ independent sets $A_1, \dots , A_m , B_1, \dots , B_{ \tau}$.  Making each one of these independent sets a color class gives (\ref{eq:general ub}).  
Let $A_1 = M_1 \times R_2 \times R_3 \times \cdots \times R_m$.  For $2 \leq i \leq m$, let 
\[
A_i = ( R_1 \times \cdots \times 
R_{i-1} \times M_i \times R_{i+1} \times \cdots \times R_m ) \backslash (A_1 \cup \dots \cup A_{i-1} ).
\]
Note that every $x\in A_i$ has its $i$-th coordinate in $M_i$, thus such an $x$ is a non-unit in $R$. Further, every non-unit in $R$ has a coordinate in $M_j$ for some minimal $j$ and must be an element of $A_j$. It follows that the union $\cup_{i=1}^m A_i$ is the set of non-units in $R$.
Let $i \in \{1,2, \dots ,m \}$.  
For any distinct $x,y \in A_i$, the $i$-th coordinate of $x+y$ is in $M_i$ so that $x+y$ is not a unit.  Therefore, each $A_i$ is an independent set. 

On the units $R^*$, define a relation $\sim$ by the rule $x \sim y$ if and only if there are integers $i_1 , \dots , i_m \in \{0,1 \}$ such that 
\[
(x_1 , \dots , x_m  ) = ( (-1)^{i_1} y_1 , \dots , (-1)^{i_m } y_m ).  
\]
It can be checked that $\sim$ is an equivalence relation.  Suppose that $B$ is an equivalence class and 
$(x_1, \dots , x_m) \in B$.  If $|B| < 2^m$,     
then we must have 
\[
( (-1)^{i_1} x_1 , \dots , (-1)^{i_m} x_m) ) 
= 
( (-1)^{j_1} x_1 , \dots , (-1)^{j_m} x_m)
\]
for two distinct choices of exponents $(i_1, \dots , i_m)$ and $(j_1 , \dots , j_m)$.  Without loss of generality, assume $i_1 = 0$ and $j_1 = 1$.  Then $x_1 = -x_1$ so that $x_1 + x_1 = 0_{R_1}$.  Since $x_1 \notin M_1$, we know that $x_1 \neq 0_{R_1}$ and so $x_1$ has additive order 2 in $R_1$.  This contradicts the fact that $R$ has an odd number of elements.  We conclude that $|B| =2^m$ and all equivalence classes have $2^m$ elements.  
Let $B_1 , \dots , B_{ \tau}$ be the equivalence classes under the relation $\sim$ on $R^*$.  We will complete the proof by showing that each $B_i$ is an independent set.
 Let $x$ and $y$ be distinct vertices in $B_i$.  
Then 
\[
(x_1 , \dots , x_m) = ( (-1)^{i_1} y_1 , \dots , 
(-1)^{i_m} y_m )
\]
for some $i_1,\dots , i_m \in \{0,1 \}$ and at least one $i_j$ is not 0. Without loss of generality, assume $i_1 = 1$.  Then the first coordinate of $x+y$ is $0_{R_1}$ and so $x+y$ is not a unit.  Therefore, $B_i$ is an independent set.  Combining the $B_i$'s with the $A_i$'s gives a proper coloring of $\mathcal{U} (R)$ with 
\[
m + \tau = m + \frac{ |R^* | }{2^m } 
\]
colors.  The proof of (\ref{eq:general ub}) is completed by noting that  
\[
|R^*| = |R_1^*|  \cdots | R_m^*|
= (|R_1| - |M_1| )  \cdots ( |R_m| - |M_m| )
=
\prod_{i=1}^m ( |R_i| - |M_i| ).
\]

Now we prove the lower bound
\begin{equation}\label{eq:general lb}
m +  \prod_{i=1}^m \left( \frac{ |R_i| - |M_i | } {2} \right)  \leq \omega ( \mathcal{U} ( R_1 \times \cdots \times R_m))
\end{equation}
by demonstrating a clique of size $m+|R^*|/2^m$. Fix an $i$ with $1 \leq i \leq m$.  
Let $\pi_i : R_i \rightarrow R_i / M_i$ be the quotient map $\pi_i (r) = r + M_i$.  Since $M_i$ is maximal, 
$R_i/M_i$ is a finite division ring and thus, is a field (or recall $R$ is commutative).  Since $|R|$ is odd, we have the $|R_i / M_i | = q_i$ for some odd prime power $q_i$.  Let $\hat{S}_i$ be a set of $\frac{q_i - 1}{2}$ non-zero elements of $R_i / M_i$ such that
\begin{enumerate}
    \item $x + M_i \in \hat{S}_i$ if and only if 
$-(x+M_i) \notin \hat{S}_i$, and
\item $-(1_{R_i} + M_i) \notin \hat{S}_i$.  
\end{enumerate}
Such a set exists since $(R_i / M_i,+)$ is an abelian group of odd order.  Let $S_i = \{x \in R_i : \pi_i ( x ) \in \hat{S}_i \}$.  

\medskip
\noindent
\textit{Claim 1:} The sum of any two distinct elements of $S_i$ is a unit in $R_i$.
\begin{proof}[Proof of Claim 1]
    Let $x,y $ be distinct elements in $S_i$.  Then $x+y$ is a unit in $R_i$ if and only if $x+y \notin M_i$.  Now if $x+y \in M_i$, then $\pi_i (x+y) = 0 + M_i$ which implies 
    \[ 
    (x+M_i ) + (y + M_i) = 0 + M_i.
    \]
    This is a contradiction since $ \hat{S}_i$ was chosen so that $x + M_i \in  \hat{S}_i$ if and only if $-(x+M_i) \notin \hat{S}_i$.  Thus, $x+y \notin M_i$.   
\end{proof}

\medskip
\noindent
\textit{Claim 2:} $|S_i| =  \left( \frac{ |R_i| - |M_i| }{2} \right)$.
\begin{proof}[Proof of Claim 2]
From the map $\pi_i$, we have $|S_i| = |M_i| ( \frac{q_i - 1}{2} )$.  Since $|R_i / M_i | = q_i$, we obtain by substitution
\[
|S_i| = \left(  \frac{ |R_i| / |M_i|  - 1 }{2} \right) |M_i| = \left( \frac{ |R_i| - |M_i| }{2} \right).  
\]
\end{proof}
By Claims 1 and 2, the set $S_1 \times \cdots \times S_m$ forms a clique of size $\prod_{i=1}^m |S_i|$.  We will define a set $T$, with $m$ elements, such that 
$T \cup ( S_1 \times \cdots \times S_m )$ is a clique.
Let
\[
T = \{ ( 0,1,1, \dots , 1,1) , (1,0,1, \dots , 1 , 1) , 
(1,1,0, \dots , 1, 1) , \dots , (1,1,1, \dots , 0,1) , 
(1,1,1, \dots , 1,0) \}.  
\]
Since $|R|$ is odd, the sum of any two elements in $T$ is a unit because $1_{R_i}$ and $1_{R_i} + 1_{R_i}$ are units for every $i$.  The next step is to show that every vertex in $T$ is adjacent to every vertex in 
$S_1 \times \cdots \times S_m$.  Without loss of generality, consider the element 
$(0,1, \dots , 1)$ in $T$.   
Let $(s_1 , s_2 , \dots , s_m)$ be in $S_1 \times \cdots \times S_m$.  If $1_{R_i} + s_i $ is not a unit in $R_i$, then $1_{R_i} + s_i \in M_i$.  Therefore, 
$-1_{R_i} + M_i = s_i + M_i$ which implies 
$- ( 1_{R_i} + M_i ) \in \hat{S}_i$, a contradiction.
The conclusion is that $1_{R_i} +s_i$ is a unit in $R_i$ for $2 \leq i \leq m$.  Since $s_1$ is a unit in $R_1$, we have that $(s_1,s_2, \dots ,s_m) + (0,1, \dots , 1)$ is a unit in $R_1 \times \cdots \times  R_m$.  This shows that $T \cup (S_1 \times \cdots \times S_m)$ is a clique.  It has size 
$$m + \prod_{i=1}^m |S_i| = 
m + \prod_{i=1}^m \left( \frac{ |R_i |  - |M_i| }{2} \right)$$
which completes the proof of Theorem \ref{theorem:general clique and chromatic}. \end{proof}


\section{Proof of Theorem \ref{thm:achromatic theorem 1}}\label{section:achromatic}

\begin{proof}[Proof of Theorem \ref{thm:achromatic theorem 1}]
 We work with the graph $\mathcal{U}( \mathbb{Z}_p \times \mathbb{Z}_q )$ which is isomorphic to   
 $\mathcal{U}( \mathbb{Z}_{pq} )$.  Let $3 \leq p < q$ be primes.  We will define an achromatic coloring of 
$\mathcal{U} ( \mathbb{Z}_p \times \mathbb{Z}_q )$
by defining a partition of the vertex set into $\frac{pq+1}{2}$ color 
classes.  Each class, with the exception of one, will contain two vertices.  

Let $C_s = \{ (0 , \frac{q-1}{2} ) \}$, 
$E_1 = \{ ( 0,0) , ( p - 1 , 0 ) \}$, $E_2 = \{ (1,0) , ( p - 1 , \frac{q+1}{2}) \}$, $E_3 = \{ (0,1 ), (p-1,q-1) \}$, $E_4 = \{ (1,1) , (0,q-1 ) \}$, $E_5 = \{ (2,1) , (1,q-1 ) \}$, and $E_6 = \{ ( p - 1 , 1) , ( p -2 , q- 1 ) \}$.  Next, let 
\begin{center}
	$M_z = \{ ( z + 2 , 1 ) , ( \frac{p-1}{2} + z , q - 1 ) \}$ for $1 \leq z \leq \frac{p-5}{2}$.
\end{center}
The set $M_z$ is empty when $p \in \{3, 5 \}$.  
For $1 \leq t \leq \frac{p-3}{2}$, let 
\begin{center}
	$D_t = \{ ( p - t , 0 ) , ( t , q- 1) \}$ ~~and ~~$D_{t}' = \{ ( t+1 , 0 ) , ( p - t - 1 , 1 ) \}$.  
\end{center}
The sets $D_t$ and $D_{t}'$ are empty when $p  = 3$.  
For $1 \leq i \leq p - 1$, let 
\begin{center}
    $C_{i, \frac{q-1}{2} } = \{  (i, \frac{q-1}{2} ) , (  i - 1  , \frac{q+1}{2} )\}$.
    \end{center}
Lastly, let 
\begin{center}
	$C_{i,j} = \{ ( i,j) , (i-1 , q - j) \}$ for $0 \leq i \leq  p - 1$ and $2 \leq j \leq \frac{q-3}{2}$.
\end{center}
This completes the definition of the color classes.  

Certainly $C_s$ is an independent set.  Every set of the form $C_{i,j}$ or $C_{i , \frac{q-1}{2} }$ or $M_z$ is an independent set because the sum of any two vertices within one of these sets has a zero $y$-coordinate.   
Similarly, every pair of vertices in a $D_t$ or a $D_{t}'$ has a sum with a zero $x$-coordinate.  It can be checked directly that each $E_i$ is also an independent set, thus, this partition defines a proper coloring of $\mathcal{U}( \mathbb{Z}_p \times \mathbb{Z}_q )$.  It remains to show that there is an edge between every pair of color classes.  

Starting with $C_s = \{ ( 0 , \frac{q-1}{2} )\}$, we list a neighbor of $(0, \frac{q-1}{2} )$ in each color class.  This is shown in the table below.   
\begin{center}
	\begin{adjustbox}{width=\textwidth}
\begin{tabular}{|c|c|c|c|c|c|c|c|c|c|} \hline
	$E_1$ & $E_2$ & $E_3$ & $E_4$ & $E_5$ & $E_6$ & $M_z$ & $D_t$ & $D_{t}'$ & $C_{ i , \frac{q-1}{2} }$  \\
	$(p-1,0)$ & $(1,0)$ & $(p-1,q-1)$ & $(1,1)$ & $(2,1)$ & $(p-1,1)$ & 
	$(z+2,1)$ & $(p-t,0)$ & $(t+1,0)$ & $(i , \frac{q-1}{2} )$ \\ \hline
		\end{tabular}
		\end{adjustbox}
		\end{center}
The last type of color class to check is a $C_{i,j} $.  
If $i \neq 0$, vertex $(i,j)$, which is in $C_{i,j}$, is adjacent to $(0 , \frac{q-1}{2} )$. If $i = 0$, then vertex $(p-1,q-j)$ is adjacent to $( 0  , \frac{q-1}{2} )$.  Indeed, since $2 \leq j \leq \frac{q-3}{2}$, the $y$-coordinate of $(p-1,q-j)$ is in the set 
$\{ q-2, q-3, \dots , \frac{q+3}{2} \}$, which does not contain $\frac{q+1}{2}$.  We conclude that there is an edge between $C_s$ and every $C_{i,j}$.  This completes the portion of the proof that shows $C_s$ has an edge to every other color class. 

Having dealt with $C_s$, observe that all remaining color classes have two vertices.  Given a pair of vertices $\{ (z_1,t_1) , (z_2,t_2) \}$, write 
$N^c ( \{ (z_1,t_1 ) , (z_2 , t_2 ) \} )$ for the set of vertices that are not adjacent to $(z_1 ,t_1)$ and not adjacent to $(z_2,t_2)$.  If $C = \{ (z_1,t_1) , (z_2,t_2) \}$ is a color class, then we cannot have another color class $C'$ contained in $ N^c ( C ) $.  Otherwise, there would be no edge between $C$ and $C'$.   This observation will be used in the arguments to follow. In particular, if $z_1 \not\equiv z_2 ( \textup{mod}~p)$ and 
$t_1 \not\equiv t_2 ( \textup{mod}~q)$, then 
\begin{equation}\label{eq:inclusion}
N^c ( \{ ( z_1,t_1) , (z_2 ,t_2) \} ) \subseteq 
\{ (  p - z_1 , q-t_2 ) , ( p - z_2 , q - t_1) \}.
\end{equation}
For some color classes, like $E_2$, $E_3$, $E_4$, and  $D_t$, $D_t'$ for $t \in \{1,2, \dots , \frac{p-3}{2} \}$, inclusion (\ref{eq:inclusion}) may be a proper containment.  The reason is that one of the vertices $(p-z_1 , q-t_2)$ or $(p- z_2 , q-t_1)$ is one of the vertices in the  color class $\{ (z_1,t_1),(z_2,t_2) \}$.
If $z_1 \equiv z_2 ( \textup{mod}~p)$ or 
$t_1 \equiv t_2 ( \textup{mod}~q)$, then (\ref{eq:inclusion}) may not hold (this occurs for color class $E_1$).  

Consider the color classes $E_1, \dots ,E_6$.  
Starting with $E_1 =\{ ( 0,0 ) , ( p - 1 , 0) \}$, 
we have $$N^c ( \{ ( 0,0) , (p -1 , 0) \} ) = \{ ( x ,0 ) : 0 \leq x \leq p \}.$$  By inspection it can be checked that there is no color class $C$, other than $E_1$, where all vertices in $C$ have $y$-coordinate equal to 0.  Hence, $N^c ( \{ ( 0,0) , (p -1 , 0) \} ) $ does not contain a color class and so there is an edge from $E_1$ to every other color class.  

Next we consider $E_2 = \{ (1,0) , ( p - 1 , \frac{q+1}{2} ) \}$.  Observe that 
\begin{center}
$N^c ( \{ ( 1,0) , (  p - 1 , \frac{q+1}{2} ) \} ) =
\{ (p-1 , \frac{q-1}{2} ) \}$.  
\end{center}
This set does not contain a color class and so there is at least one edge between $E_2$ and each of the other color classes.  Similarly, the set
\[
N^c ( E_3 ) = N^c ( \{ (0,1 ) , ( p - 1 , q- 1) \} ) = \{ (1 ,  q - 1 ) \}
\]
does not contain a color class, hence there is at least one edge between $E_3$ and any other color class.  
The same is true for $E_4$ since 
\[
N^c (E_4) = N^c ( \{ ( 1,1) , ( 0 , q - 1) \} ) 
=
\{  (p-1,1) \}.    
\]

With $E_5 = \{ (2,1) , (1,q-1) \} $, we have 
\[
N^c ( E_5) = \{ ( p - 1 ,  q- 1) , ( p -2 , 1) \}.  
\]
Vertex $(p-1,q-1)$ is in $E_3 = \{ ( 0,1 ) , ( p - 1 , q- 1 ) \}$.  Hence, $N^c ( E_5)$ does not contain a color class.  Similarly, 
\[
N^c (E_6) = N^c ( \{ ( p-1 ,1) , (p-2,q-1) \} ) 
= \{ (1,1) , (2,q-1) \}. 
\]
Vertex $(1,1)$ is in color class $E_4  = \{ ( 1,1) , (0,q-1 ) \}$.  Thus, $N^c ( E_6)$ does not contain a color class.  

Let $0 \leq i \leq p - 1$ and $2 \leq j \leq \frac{q-3}{2}$.  Then 
\[
N^c (C_{i,j} ) = 
N^c ( \{  ( i , j ) , (i - 1 , q- j ) \} ) 
= \{ ( p - i , j ) , ( p - ( i -1 ) , q - j ) \}.
\]
Vertex $(p-i,j)$ is in color class $C_{p-i,j}$ while $(p-(i-1) , q-j )$ is in $C_{p-(i-1) , j } $.  These are distinct color classes for these values of $i$ and $j$.

For $1 \leq i \leq p - 1$, 
\begin{center}
$N^c ( C_{i , \frac{q-1}{2}  } ) 
= 
N^c ( \{ ( i , \frac{q-1}{2} ) , ( i - 1 , \frac{q+1}{2} ) \} ) 
= \{ ( p - i , \frac{q-1}{2} ) , ( p - ( i - 1 ) ,\frac{q+1}{2} ) \}$.
\end{center}
Vertex $(p-i, \frac{q-1}{2} )$ is in the color class $C_{p-i, \frac{q-1}{2} } = \{ ( p - i , \frac{q-1}{2} ) , ( p - i - 1  , \frac{q+1}{2} ) \}$.  Therefore, $N^c ( C_{i, \frac{q-1}{2} } )$ does not contain a color class. 

Let $1 \leq t \leq \frac{p-3}{2}$.  In this case, 
\[
N^c ( D_t) = N^c ( \{ ( p - t , 0 ) , ( t,  q - 1 ) \} ) 
 = 
 \{ ( t , 1 ) \}
\]
and 
\[
N^c ( D_{t}') = N^c ( \{ ( t+1 , 0 ) , ( p-t-1,1)  \} ) 
 = 
 \{ ( p-(t+1) , q -1 )  \}.
\]
Recalling that the only color class that consists of one vertex is $C_s = \{ ( 0 , \frac{q-1}{2} \}$, we may conclude that the color classes $D_t$ and $D_t$' send at least one edge to the other color classes.

Assume now that $p > 5$ so that $M_z$ is not empty.  
Let $1 \leq z \leq \frac{p-5}{2}$.  We have 
\begin{center}
$N^c ( M_z ) = 
N^c ( \{ ( z + 2 , 1 ) , ( \frac{p-1}{2} + z , q - 1 ) \} ) 
= 
\{ ( p - ( z + 2) , 1) , (  \frac{p+1}{2} -z , q - 1 ) \}$.  
\end{center}
The vertex $(p - ( z + 2),1)$ is in the color class 
\begin{center}
    $M_{ p - z - 4} = \{ ( p - ( z + 2) , 1 ) , ( \frac{p-1}{2} + p - z - 4 , q - 1 ) \}$.
    \end{center}
The assumption that $p > 5$ implies that 
the vertices $ ( \frac{p-1}{2} + p - z - 4 , q - 1 )$
and $ (  \frac{p+1}{2} -z , q - 1 ) $ are not equal.  We conclude that $N^c( M_z)$ does not contain a color class.  

To complete the proof of Theorem \ref{thm:achromatic theorem 1}, we will compute the number of color classes. 
 For $p \geq 5$, this value is
 \[ 
 1 + 6 + \frac{p-5}{2} + 2 \left( \frac{p-3}{2} \right) 
 + (p-1) 
 +
 p \left( \frac{ q-5}{2} \right) = \frac{pq+1}{2}.
\]
In the case that $p = 3$, there are no $D_t$ or $D_t'$ color classes.  Furthermore, $E_5 = E_6$ in this case.  Hence, the number of color classes when $ p  = 3$ is 
\[
1 + 5 + 2 + 3 \left( \frac{q-5}{2} \right) = \frac{3q+1}{2}.  
\]
  This completes the proof for all primes $q > p  \geq 3$.  
\end{proof}


\section{Proof of Theorem \ref{thm:achromatic theorem 2}}

We begin this section by determining $\chi_a ( \mathcal{U} ( \mathbb{Z}_{2p} )) $ as it provides motivation for a definition we give in a moment. 
 The argument is based on the maximum degree of the graph $\mathcal{U} ( \mathbb{Z}_{2p} )$.  We will also use a result of Harary, Hedetniemi and Prins (Theorem 3 in \cite{Harary}) which states that for any graph $G$ and integer $t$ with $\chi (G) \leq t \leq \chi_a (G)$, $G$ has an achromatic coloring with $t$ colors.  

\begin{proof}[Proof of $\chi_a( \mathcal{U} ( \mathbb{Z}_{2p} ) )  ) = p $] Aiming for a contradiction, suppose that 
$\chi_a ( \mathcal{U} ( \mathbb{Z}_{2p} ) ) > p$, so that 
$\mathcal{U}( \mathbb{Z}_{2p} )$ has an
achromatic coloring with $p+1$ colors.  
Since this graph has $2p$ vertices, the average size of a color class is $\frac{2p}{p+1} < 2$.  
Therefore, there is a color class that contains only one vertex.  Let $v$ be such a vertex.  The degree of $v$ is at most $\phi (2p) = p-1$, so 
that $v$ cannot have at least one neighbor in each of the remaining $p$ color classes.  This is a contradiction, so we must have 
\[
\chi_a( \mathcal{U} ( \mathbb{Z}_{2p} ) )  ) = p 
\]
for any odd prime $p$ (the lower bound 
$\chi_a ( \mathcal{U} ( \mathbb{Z}_{2p} ) ) \geq p$ is proved in \cite{Momrit}).  \end{proof}

Now we turn our attention to the main part of this section which is the proof of Theorem \ref{thm:achromatic theorem 2}.  Let us start with a definition.  

\begin{defn} Given an achromatic coloring of a graph, a vertex is \textbf{special} for that coloring if it is the only vertex in its color class.
\end{defn}
Observe that in any achromatic coloring of a graph, the set of special vertices form a clique.  

\begin{proof}[Proof of Theorem \ref{thm:achromatic theorem 2}] 
Aiming for a contradiction, suppose we have an achromatic coloring of $\mathcal{U} ( \mathbb{Z}_3 \times \mathbb{Z}_q ) $ with 
$\frac{3q+1}{2} + 1$ colors.  Call this coloring $f$.  
If there are at most two special vertices, then since 
 $\mathcal{U} ( \mathbb{Z}_3 \times \mathbb{Z}_q ) $ has $3q$ vertices, the number of color classes is at most 
\[
2 + \frac{3q-2}{2} = \frac{3q+2}{2} < \frac{3q+3}{2}.
\]
This is a contradiction because $f$ uses $\frac{3q+3}{2}$ colors.  
Thus, there are at least three special vertices.  Let
\[
\mathcal{S} = \{ (x_1,y_1) , (x_2,y_2)  ,  \dots , (x_k,y_k) \}
\]
be the set of all special vertices under $f$.  
Since the vertices in $\mathcal{S}$ form a clique, 
\begin{center}
    $x_i + x_j \not\equiv 0 ( \textup{mod}~3)$ and $y_i + y_j \not\equiv 0 ( \textup{mod}~q)$ for $1 \leq i < j \leq k$.
    \end{center}
    In particular, at most one $x_i$ and one $y_j$ are 0.  The function that maps the vertex $(z,t)$ to $(-z,t)$ is a graph automorphism of $\mathcal{U} ( \mathbb{Z}_3 \times \mathbb{Z}_q ) $.  Thus, we may assume that $x_i = 1$ for at least one $i$.  Again, using the fact that $\mathcal{S}$ is a clique, we may now say that $x_i \in \{0,1 \}$ for $1 \leq i \leq k$ and $x_i = 0$ for at most one $i$.

    For $i \in \{1,2, \dots , k \}$, let $C( (0,-y_i) )$ be the color class that contains $(0,-y_i)$.  
Let $(x_i ,y_i)$ be a special vertex for which $x_i =1$.  At least two such vertices exist.  
In this case, $(1,y_i) \neq (0,-y_i)$, and so
the color class $C( (0,-y_i ))$ must contain a neighbor of $(1,y_i)$. Let $(u_i ,v_i)$ be a neighbor of $(1 , y_i)$ that is in $C( (0 , - y_i))$.  
Since $(u_i,v_i) + ( 1 , y_i)$ is a unit, $u_i \in \{0,1 \}$ and also, $v_i \not\equiv -y_i ( \textup{mod}~q)$.  
If $u_i = 1$, then for $(1,v_i) + ( 0 , - y_i)$ to be a non-unit, we require $v_i \equiv y_i ( \textup{mod}~q)$.  This is a contradiction since $(u_i,v_i)$ cannot be the same vertex as $  (1,y_i)$.  Therefore, $u_i = 0$.  For $(0,v_i)  + ( 1 , y_i)$ to be a unit, we require
$v_i \not\equiv -y_i ( \textup{mod}~q)$.  
At this point, we have shown that the color class $C((0,-y_i))$ contains $(0,-y_i)$ and $(0,v_i)$ for some $v_i \not\equiv -y_i ( \textup{mod}~q)$ in the case that our special vertex $(x_i,y_i)$ has $x_i = 1$.  Let us now prove a property about all vertices in this color class.  

 Suppose that $C((0,-y_i) )$ contains some vertex $(z,t)$ with $z \in \{1,2 \}$.  At least one of the sums 
 $(z,t) + (0,-y_i)= (z , t -y_i)$ or 
 $(z,t) + (0,v_i) = (z , t + v_i)$ will be a unit because not both 
 $t-y_i$ and $t + v_i$ can be congruent to zero mod $q$ for one value of $t$ (recall $v_i \neq - y_i$).  We conclude that all vertices in $C((0,-y_i))$ are of the form $(0,t)$ for some $t \in \mathbb{Z}_q$.   In other words,  
 \begin{equation}\label{new eq 5}
 C( (0,-y_i )) \subseteq \{ 0 \} \times \mathbb{Z}_q.
 \end{equation}
 
If $x_i = 1$ and $x_j = 1$ for some $i,j \in \{1,2, \dots , k \}$, then the sum of any two vertices 
in $C( (0, - y_i ) ) \cup C ( ( 0 , - y_j ))$ will have a 0 $x$-coordinate, and that pair of vertices will not be adjacent.  In order for there to be an edge between any two color classes, all of the $C((0,-y_i))$'s for which (\ref{new eq 5}) holds must be the same class.  Let $\mathcal{C}_0$ denote this color class.  Since $k \geq 3$, there are at least two distinct $i$'s for which $x_i = 1$.  Thus, $\mathcal{C}_0$ is not empty.  This implies that there can be no special vertex $(x_i,y_i)$ with $x_i = 0$ since such a vertex does not have a neighbor in $\mathcal{C}_0$.  Hence,
\[
\mathcal{S} = \{ ( 1,y_1) , (1,y_2) , \dots , (1,y_k) \}.
\]
Let us now count the number of colors that $f$ uses on the vertices in 
$\mathcal{C}_o \cup \mathcal{S}$ and its complement.  

 The coloring $f$ uses $1+k$ colors on the color class $\mathcal{C}_0$ and the special vertices $(1,y_1), \dots , (1,y_k)$.  This leaves $\frac{3q+1}{2} + 1 - ( 1 +k )$ colors for the remaining $3q - k - |\mathcal{C}_0|$ vertices, none of which is special.  
 The color classes containing these vertices will all have more than one vertex.  Hence, 
 the number of vertices contained in the union of these $\frac{3q+1}{2} - k$ color classes is at least
\[
2 \left( \frac{3q+1}{2}  - k \right) = 3q + 1 - 2k.
\]
Thus,
\[
3q - k - |\mathcal{C}_0 | \geq 3q + 1 - 2k.
\]
Solving for $|\mathcal{C}_0 |$ gives $k-1 \geq |\mathcal{C}_0 |$.  This contradicts the fact that the color class $\mathcal{C}_0 $ contains the $k$ distinct vertices $(0,-y_1), \dots , (0,-y_k)$.  
We conclude that no such coloring $f$ exists which  completes the proof.  
\end{proof} 

\section{Concluding Remarks}\label{se:concluding remarks}

One of our reasons for studying the achromatic number of $\mathcal{U}( \mathbb{Z}_n)$ is attempting to solve a problem of Ahlswede, Bezrukov, Blokhuis, Metsch, and Moorhouse \cite{Ahlswede}.  Given a graph $G$, let $\psi (G)$ be the maximum number of colors in a, not necessarily proper, vertex coloring of $G$ such that there is at least one edge between every pair of distinct color classes.  The value $\psi (G)$ is often called the \emph{pseudo-achromatic number} of $G$.   
For any graph $G$, one has $\chi (G) \leq \chi_a (G) \leq \psi (G)$.  Also, since there must be an edge between any two color classes, 
\begin{center}
    $\dbinom{\psi (G) }{2} \leq e(G)$.
   \end{center}
   The same inequality holds with $\chi (G)$ or $\chi_a(G)$ in place of $\psi (G)$.  It implies 
   \[
   \psi (G) \leq \sqrt{ 2 e(G) + 1/4} + 1/2.
   \]
   An interesting question is given a graph $G$, how close are $\chi_a(G)$ and $\psi (G)$ to 
   $\sqrt{ 2 e(G) + 1/4} + 1/2$.
   Let $\{ G_n \}_{n=1}^{ \infty} $ be a sequence of graphs such that $|V(G_n) | \rightarrow \infty$ as $n \rightarrow \infty$.  Ahlswede et al. \cite{Ahlswede} asked for an explicit sequence of graphs $G_n$ such that 
\begin{equation}\label{A}
\liminf_{n \rightarrow \infty} \frac{ \psi (G_n) } { \sqrt{ 2 e(G_n) } } = 0 .
\end{equation}
To our knowledge, no explicit sequence of graphs $G_n$ that satisfy (\ref{A}) is known.  Roichman \cite{Roichman} called a sequence $G_n$ \emph{almost optimally complete} if 
   \[
   \lim_{n \rightarrow \infty} \frac{ \chi_a ( G_n ) }{ \sqrt{ 2 e(G_n) } }  = 1.  
   \]
It is known that the sequence of paths $\{P_1 , P_2 , \dots \}$ and the sequence of cycles $\{C_3 , C_4 , \dots \}$ are almost optimally complete.  More generally, any sequence of graphs of bounded degree is almost optimally complete \cite{Edwards}.  
By Theorem \ref{thm:achromatic theorem 2},
$\chi_a ( \mathcal{U} ( \mathbb{Z}_{3q} ) )  = \frac{3q+1}{2}$ for all primes $q \geq 5$.  There is roughly a factor of $\frac{8}{3}$ between 
$\binom{ \chi_a ( \mathcal{U} ( \mathbb{Z}_{3q} ) ) }{2} = \frac{9q^2-1}{8}$ and 
$e ( \mathcal{U} ( \mathbb{Z}_{3q} ) ) = \frac{1}{2}(3q-1)\cdot 2 (q-1) =(3q-1)(q-1)$.  
Therefore, the sequence $\{ \mathcal{U} ( \mathbb{Z}_{3q} ) : q = 5,7,11,\dots \}$ is not almost optimally complete.  

In light of these questions of Ahlswede et al.\ \cite{Ahlswede} and Roichman \cite{Roichman}, an interesting sequence to investigate is $\{ \mathcal{U} ( \mathbb{Z}_{p_1p_2 \cdots p _m} ) \}_{m=1}^{ \infty}$
where $3 = p_1 < p_2 < \dots < p_m  < \dots$ is the sequence of odd primes.  We can lower bound the number of edges in $\mathcal{U} ( \mathbb{Z}_{p_1p_2 \dots p _m} )$  using results in number theory.    
Writing $p_N$ for the $N$-th prime, the work of Rosser and Shoenfeld \cite{Rosser} gives $p_N \leq 2N \log N$ for all $N \geq 6$, and
$\prod_{p \leq x } \left( 1 - \frac{1}{p} \right) 
> 
\frac{ \textup{exp}(  - \gamma ) }{ \log x } \left( 1 - \frac{1}{ 2 \log^2 x } \right)$ for $x \geq 285$.  Here the product is over all prime numbers that are less than or equal to $x$, $\gamma$ is the Euler-Mascheroni constant, and $\textup{exp}(t):=e^t$.  Using these inequalities we have for large enough $m$, 
\[
\prod_{i=1}^m \left( 1 - \frac{1}{p_i} \right) 
=
2 \prod_{p \leq p_m} \left( 1 - \frac{1}{p} \right)
>
\frac{2 \textup{exp}( - \gamma ) }{ \log (p_m) } 
\left( 1 - \frac{1}{2 \log^2 ( p_m) } \right) 
\geq 
\frac{ 2 \textup{exp} (  - \gamma ) }{ \log (2m \log m) }.
\]
Hence, 
\begin{eqnarray*}
e( \mathcal{U} ( \mathbb{Z}_{p_1p_2 \dots p _m} ) ) 
&=& 
\frac{1}{2} ( p_1 p_2 \cdots p_m - 1 ) \phi (p_1 p_2 \cdots p_m)  \\
 &= &  
\frac{1}{2} (p_1 p_2 \cdots p_m - 1 ) p_1 p_2 \cdots p_m
\prod_{i=1}^m \left( 1 - \frac{1}{p_i} \right)
\\
& > & 
\frac{1}{2}(p_1 p_2 \cdots p_m -1 )^2 \frac{ \textup{exp}(  - \gamma ) }{ \log ( 2m \log m ) }
\end{eqnarray*}
which gives 
\[
\frac{
\psi (  \mathcal{U} ( \mathbb{Z}_{p_1p_2 \dots p _m} ) ) }
{ \sqrt{ 2 e( \mathcal{U} ( \mathbb{Z}_{p_1p_2 \dots p _m} ) ) }}
\leq
\frac{\psi (  \mathcal{U} ( \mathbb{Z}_{p_1p_2 \dots p _m} ) ) 
 \textup{exp}(   \gamma /2 )  \sqrt{ \log ( 2m \log m ) } }{  p_1 p_2 \cdots p_m - 1  }. 
\]
This leads us to ask whether or not 
$\psi (  \mathcal{U} ( \mathbb{Z}_{p_1p_2 \dots p _m} ) )$
is $o ( \frac{p_1 p_2 \cdots p_m  }{ \sqrt{ \log (m \log m) } } )$?  Since a pseudo-achromatic coloring need not be proper, proving a good upper bound on 
$\psi (  \mathcal{U} ( \mathbb{Z}_{p_1p_2 \dots p _m} ))$ 
could be more difficult than proving a good upper bound 
on $\chi_a (  \mathcal{U} ( \mathbb{Z}_{p_1p_2 \dots p _m} ))$.  In this paper all colorings considered were proper and this assumption was used in the proof of Theorem \ref{thm:achromatic theorem 2}.  



\begin{thebibliography}{99}

\bibitem{Ahlswede}
R.\ Ahlswede, S.\ L.\ Bezrukov, 
A.\ Blokhuis, K.\ Metsch, G.\ E.\ Moorhouse, 
Partitioning the $n$-cube into sets with mutual distance 1,
{\em Appl.\ Math.\ Lett}. Vol.\ 6, No.\ 4 (1993), 17--19.  

\bibitem{Akhtar}
R.\ Akhtar, M.\ Boggess, 
T.\ Jackson-Henderson, I.\ Jim\'{e}nez, 
R.\ Karpman, A.\ Kinzel, D.\ Pritikim, 
On the unitary Cayley graph of a finite ring, 
{\em Electron.\ J.\ Combin.}, {\bf 16} (2009), R117.  

\bibitem{Anderson}
D.\ Anderson, T.\ Asir, A.\ Badawi, T.\ Tamizh Chelvam,
{\em Graphs from Rings}, 
Nov 2021, Springer Nature Switzerland AG 2021.  



 \bibitem{Bahrami}
 A.\ Bahrami, R.\ Jahani-Nezhad,
 Unit and unitary Cayley graphs for the ring of Gaussian integers modulo $n$, 
 {\em Quasigraphs and Related Systems}, 
 25 (2017), 189--200.  

\bibitem{Chin} 
A.\ Y.\ M.\ Chin, H.\ R.\ Maimani, M.\ R.\ Pournaki, M.\ Sivagami, T.\ Tamizh Chelvam,
Unitary Cayley graphs whose Roman domination numbers are at most four
{\em AKCE International J.\ of Graphs and Combinatorics}, {\bf 19(1)} (2022), 36--40.   

\bibitem{Edwards}
N.\ Cairnie, K.\ Edwards, 
Some results on the achromatic number
{\em J.\ Graph Theory}, {\bf 26}| (1997) 129--136.

\bibitem{EdwardsSurvey}
K.\ Edwards, 
The harmonious chromatic number and the achromatic number,
{\em Surveys in Combinatorics, 1997}, London Mathematical Society Lecture Note Series. Cambridge University Press (1997), 13--48.  
  

\bibitem{Harary2}
F.\ Harary, S.\ Hedetniemi,
The achromatic number of a graph,
{\em J.\ Combin.\ Theory  Ser.\ B}, {\bf 8} (1970), 154--161.  

\bibitem{Harary}
F.\ Harary, S.\ Hedetniemi, G.\ Prins,
An interpolation theorem for graphical homomorphisms, Portugal.\ Math.\ 26 (1967) 453--462.  

\bibitem{Klotz}
W.\ Klotz, T.\ Sander
Some properties of unitary Cayley graphs,
{\em Electron.\ J.\ Combin.}, {\bf 14}, R45 (2007).  

\bibitem{Liu}
X.\ Liu, S.\ Zhou,
Eigenvalues of Cayley graphs, 
{\em Electron.\ J.\ Combin.}, {\bf 29(2)}, P2.9 (2022).  
  
\bibitem{Momrit}
P.\ Momrit, C.\ Promsakon, 
The achromatic numbers of unitary addition Cayley graphs, 
{\em Proceedings of AMM 2017}, The 22nd Annual Meeting in Mathematics.
		
\bibitem{Palanivel}
N.\ Palanivel, A.\ V.\ Chithra, 
Some structural properties of unitary addition Cayley graphs
{\em International J.\ of Computer Applications}, 121 No.\ 17, 2015.   2366--2377.  
		
\bibitem{Promsakon}
C.\ Promsakon, 
Colorability of unitary addition Cayley graphs, 
{\em Far East Journal of Mathematical Sciences}, Vol.\ 100, No.\ 2, 2016, 227--242. 

\bibitem{Roichman}
		Y.\ Roichman,
		On the achromatic number of hypercubes,
		{\em J.\ Combin.\ Theory Ser.\ B}, 79 (2000), 177--182.  


\bibitem{Rosser}
J.\ B.\ Rosser, L.\ Schoenfeld,
Approximate formulas for some functions of prime numbers, 
{\em Illinois J.\ Math}.\ {\bf 6} (1962), 64--94.  

  
\bibitem{Sinha}
D.\ Sinha, P.\ Garg, A.\ Singh, 
Some properties of unitary addition Cayley graphs, 
{\em Notes on Number Theory and Discrete Mathematics},
Vol.\ 17 (2011) 3, 49--59.  		
		
	\end{thebibliography}
\end{document}